\documentclass[11pt]{article}

\usepackage{amscd,amsmath, amssymb, fancyhdr, epsfig, color, tocloft, mathtools}
\usepackage{dutchcal}

\usepackage[backref=page]{hyperref}
\renewcommand*{\backrefalt}[4]{%
	\ifcase #1 (Not cited.)%
	\or        (Cited on page~#2.)%
	\else      (Cited on pages~#2.)%
	\fi}

\hypersetup{
	colorlinks   = true,
	citecolor    = magenta,
	linkcolor    = blue,
	urlcolor     = magenta	
}

\newcommand{\version}{version 2.0,\ \ January 8, 2021}
\makeatletter
\def\x@arrow{\DOTSB\Relbar}
\def\xlongequalsignfill@{\arrowfill@\x@arrow\Relbar\x@arrow}
\providecommand{\xlongequal}[2][]{%
	\ext@arrow 0099\xlongequalsignfill@{#1}{#2}}
\def\xlongrightarrowfill@{\arrowfill@\relbar\relbar\longrightarrow}
\newcommand{\xlongrightarrow}[2][]{%
	\ext@arrow 0099\xlongrightarrowfill@{#1}{#2}}
\makeatother

\numberwithin{equation}{section}

\def\eqref#1{(\ref{#1})}
\newcommand{\goth}{\mathfrak}

\newcommand{\Z}{{\mathbb Z}}
\newcommand{\C}{{\mathbb C}}
\newcommand{\R}{{\mathbb R}}
\newcommand{\Q}{{\mathbb Q}}

\newcommand{\6}{\partial}
\def\1{\sqrt{-1}\:}
\newcommand{\restrict}[1]{{\left|_{{\phantom{|}\!\!}_{#1}}\right.}}
\newcommand{\cntrct}                % contraction with a vector field
{\hspace{2pt}\raisebox{1pt}{\text{$\lrcorner$}}\hspace{2pt}}
\newcommand{\arrow}{{\:\longrightarrow\:}}

\newcommand{\calo}{{\cal O}}
\newcommand{\cac}{{\cal C}}

\renewcommand{\bar}{\overline}
\renewcommand{\phi}{\varphi}
\renewcommand{\epsilon}{\varepsilon}
\renewcommand{\geq}{\geqslant}
\renewcommand{\leq}{\leqslant}

\newcommand{\End}{\operatorname{End}}

\newcommand{\Id}{\operatorname{Id}}
\newcommand{\Teich}{\operatorname{Teich}}
\newcommand{\Comp}{\operatorname{Comp}}
\newcommand{\Gr}{\operatorname{Gr}}

\newcommand{\const}{\operatorname{\text{\sf const}}}

\newcommand{\Aut}{\operatorname{Aut}}

\newcommand{\Diff}{\operatorname{Diff}}

\newcommand{\Lie}{\operatorname{Lie}}

\newcommand{\GL}{\operatorname{GL}}

\renewcommand{\Im}{\operatorname{Im}}

\newcommand{\ie}{{\em id est }}

\newcounter{Mycounter}[section]
\newcounter{lemma}[section]
\renewcommand{\thelemma}{{Lemma \thesection.\arabic{lemma}}}
\newcommand{\lemma}{%
	\setcounter{lemma}{\value{Mycounter}}
	\refstepcounter{lemma}
	\stepcounter{Mycounter}
	{\noindent \bf \thelemma:\ }}

\newcounter{claim}[section]
\renewcommand{\theclaim}{{Claim \thesection.\arabic{claim}}}
\newcommand{\claim}{%
	\setcounter{claim}{\value{Mycounter}}
	\refstepcounter{claim}
	\stepcounter{Mycounter}
	{\noindent \bf \theclaim:\ }}

\newcounter{sublemma}[section]
\newcounter{corollary}[section]
\renewcommand{\thecorollary}{{Corollary \thesection.\arabic{corollary}}}
\newcommand{\corollary}{%
	\setcounter{corollary}{\value{Mycounter}}
	\refstepcounter{corollary}
	\stepcounter{Mycounter}
	{\noindent \bf \thecorollary:\ }}

\newcounter{theorem}[section]
\renewcommand{\thetheorem}{{Theorem \thesection.\arabic{theorem}}}
\newcommand{\theorem}{%
	\setcounter{theorem}{\value{Mycounter}}
	\refstepcounter{theorem}
	\stepcounter{Mycounter}
	{\noindent \bf \thetheorem:\ }}

\newcounter{conjecture}[section]
\newcounter{proposition}[section]
\renewcommand{\theproposition} {{Proposition \thesection.\arabic{proposition}}}
\newcommand{\proposition}{%
	\setcounter{proposition}{\value{Mycounter}}
	\refstepcounter{proposition}
	\stepcounter{Mycounter}
	{\noindent \bf \theproposition:\ }}

\newcounter{definition}[section]
\renewcommand{\thedefinition} {{Definition~\thesection.\arabic{definition}}}
\newcommand{\definition}{%
	\setcounter{definition}{\value{Mycounter}}
	\refstepcounter{definition}
	\stepcounter{Mycounter}
	{\noindent \bf \thedefinition:\ }}

\newcounter{example}[section]
\renewcommand{\theexample}{{Example \thesection.\arabic{example}}}
\newcommand{\example}{%
	\setcounter{example}{\value{Mycounter}}
	\refstepcounter{example}
	\stepcounter{Mycounter}
	{\noindent \bf \theexample:\ }}

\newcounter{remark}[section]
\renewcommand{\theremark}{{Remark \thesection.\arabic{remark}}}
\newcommand{\remark}{%
	\setcounter{remark}{\value{Mycounter}}
	\refstepcounter{remark}
	\stepcounter{Mycounter}
	{\noindent \bf \theremark:\ }}

\newcounter{problem}[section]
\newcounter{question}[section]
\makeatletter

\@addtoreset{equation}{section}
\@addtoreset{footnote}{section}
\makeatother

\def\blacksquare{\hbox{\vrule width 5pt height 5pt depth 0pt}}
\def\endproof{\blacksquare}

\newcommand{\proof}{{\bf Proof: \ }}
\newcommand{\pstep}{{\bf Proof. Step 1: \ }}

\begin{document}
	
	\begin{center}
		{\Large\bf  Lee classes on LCK manifolds with potential}\\[5mm]
		%%%%%%%%%%%%%%%%%%%%%%%%%%%%%%%%%%%%%%%%%%%%%%%%%%%%%%%%%%%%
		{\large
Liviu Ornea\footnote{Liviu Ornea is  partially supported by Romanian Ministry of Education and Research, Program PN-III, Project number PN-III-P4-ID-PCE-2020-0025, Contract  30/04.02.2021},  
Misha Verbitsky\footnote{Misha Verbitsky is partially supported by
by the HSE University Basic Research Program, FAPERJ E-26/202.912/2018 
and CNPq - Process 313608/2017-2.\\[1mm]
\noindent{\bf Keywords:} Locally conformally K\"ahler, LCK potential, Vaisman manifold, deformation, Lee form, Lee class, Hodge decomposition, algebraic cone, Teichm\"uller space.
				
\noindent {\bf 2010 Mathematics Subject Classification:} {53C55, 32G05.}
			}\\[4mm]
			
		}
		
	\end{center}

	{\small
		\hspace{0.15\linewidth}
		\begin{minipage}[t]{0.7\linewidth}
			{\bf Abstract} \\ 
An LCK manifold is a complex manifold $(M,I)$ equipped
with a Hermitian form $\omega$ and a closed 1-form
$\theta$, called the Lee form, such that
$d\omega=\theta\wedge\omega$. An LCK manifold with
potential is an LCK manifold with a positive 
K\"ahler potential on its cover, such that the deck 
group multiplies the K\"ahler potential by a constant.
A Lee class of an LCK manifold is the cohomology class of the Lee form. 
We determine the set of Lee classes on LCK manifolds 
admitting an LCK structure with potential,
showing that it is an open half-space in $H^1(M,\R)$. 
For Vaisman manifolds, this theorem was proven
in 1994 by Tsukada; we give a new self-contained proof of his 
result.
\end{minipage}
	}
	%%%%%%%%%%%%%%%%%%%%%%%%%%%%%%%%%%%%%%%%%%%%%%%%

	\tableofcontents
	
	%%%%%%%%%%%%%%%%%%%%%%%%%%%%%%%%%%%%%%%%%%%%%%%%%%%
	
	\section{Introduction}
	\label{_Intro_Section_}

%%%%%%%%%%%%%%%%%%%%%%%%%%%%%%%%%%%%%%%%%%%%%%%%

A locally conformally K\"ahler (LCK)  manifold is an
Hermitian manifold $(M,I, g)$ equipped with an atlas
$\{U_\alpha\}$ such that the restriction of $g$ to each
$U_\alpha$ is conformally equivalent to some K\"ahler
metric $g_\alpha$ defined only on $U_\alpha$, \ie
$g\restrict{U_\alpha}=e^{f_\alpha}g_\alpha$, where
$f_\alpha\in C^\infty U_\alpha$. One can see that in this
case the exterior derivatives of the conformal factors
agree on intersections: $d f_\alpha=d f_\beta$ on
$U\alpha\cap U_\beta$, thus giving rise to a global closed
1-form $\theta$, called {\em the Lee form}.
Then the Hermitian form $\omega(x,y):=g(Ix,y)$ satisfies
$d\omega=\theta \wedge \omega$ (see \cite{do} for an
introduction to the subject).

Forgetting the complex structure, one arrives at the notion of {\em locally conformally symplectic manifold} (LCS, for short): a $2n$ dimensional real manifold endowed with a non-degenerate 2-form $\omega$ and a closed 1-form $\theta$ (also called Lee form) such that $d\omega=\theta \wedge \omega$. 

Two subclasses of LCK manifolds are very important and rather well understood by now. The {\em Vaisman manifolds}, whose universal covers are K\"ahler cones over Sasakian manifolds (see Subsection \ref{_Vaisman_manifolds_}), and {\em LCK manifolds with potential}, whose universal cover admits a K\"ahler metric with global, positive and automorphic potential (see Subsection \ref{_LCK_w_potential_} for the precise definition). 

The cohomology class of the Lee form, called  {\em the Lee class}, is the first cohomological invariant one encounters when
dealing with the LCK manifolds. Let
$(M, \theta,\omega)$ be a compact LCK manifold,
and $[\theta] \in H^1(M, \R)$ its Lee class.
By Vaisman's theorem (\ref{vailcknotk}), $[\theta]=0$ if
and only if $M$ is of K\"ahler type.

For a compact  K\"ahler manifold $X$, the subset of K\"ahler classes in $H^2(X, \R)$ is the  ``K\"ahler cone'', and is
one of the most important geometric features of 
a K\"ahler manifold. Similarly, we would like to have a description of  the set 
of Lee classes on a given compact complex manifold which is known to admit LCK structures.  
It was already  shown that in this case it cannot be a cone: indeed, by A. Otiman
(\cite[Theorem 3.11]{oti2}), for an
Inoue surface of class $S^0$, the set
of Lee classes is a point. 

For LCS structures,
the set of the Lee classes is better understood, due to
Eliashberg and Murphy, who proved that 
on any almost complex manifold with $H^1(M, \Q)\neq 0$, 
for any non-zero class $\alpha\in H^1(M, \Q)$, 
there exists $C>0$ such that $C\alpha$ is the 
Lee class of an LCS structure (\cite[Theorem 1.11]{em}).

For complex surfaces with $b_1(M)=1$, the set ${\goth L}$ of 
Lee classes of LCK structures was studied by
Apostolov and Dloussky, who proved that
${\goth L}$ is either open or a point, \cite{ad2}.

For higher dimensional LCK manifolds, the first important advance in this direction 
was due to K. Tsukada, who proved that the set of Lee classes on 
Vaisman manifolds is an open half-space (\cite[Theorem 5.1]{tsuk}), using the harmonic
decomposition on Vaisman manifolds, due to T. Kashiwada, \cite{kashiwada_kodai}.

In this paper, we extend Tsukada's theorem to compact LCK
manifolds with potential of complex dimension greater than
3, using the following decomposition theorem for the first
cohomology (\ref{_LCK_pot_Hodge_decompo_Theorem_}): 
\begin{equation}\label{_H^1_decompo_Equation_}
	H^1(M, \C) = H^{1,0}(M) \oplus \overline{H^{1,0}(M)} \oplus \langle \theta \rangle
\end{equation}
where $H^{1,0}(M)\subset H^1(M, \C)$ is the space 
of all closed holomorphic 1-forms, identified with a 
subspace in cohomology by \ref{_H^1_holo_LCK_Lemma_}. Tsukada proved this
for Vaisman manifolds using the commutation formulae for Laplacians,
and the harmonic decomposition for Vaisman manifolds.

I. Vaisman conjectured that $b_1(M)$ is odd-dimensional for
any compact LCK manifold (\cite[p. 535]{va_tr}); 
this famous conjecture was disproven by
Oeljeklaus and Toma in \cite{ot}. The decomposition
\eqref{_H^1_decompo_Equation_} would imply that $b_1(M)$ is odd, hence
the counterexample of Oeljeklaus-Toma does not satisfy 
\eqref{_H^1_decompo_Equation_}. However, the natural map
\begin{equation}\label{_H^1_holom_map_Equation_}
	H^{1,0}(M) \oplus \overline{H^{1,0}(M)} \oplus \langle \theta \rangle
	\arrow H^1(M, \C)
\end{equation}
is always injective (\ref{_H^1_holo_LCK_Lemma_}).

For LCK manifolds with potential, we deduce 
\eqref{_H^1_decompo_Equation_} from a deformation argument,
by showing that an LCK manifold with potential $M_1$ obtained as a deformation
of a Vaisman manifold $M_2$ satisfies $\dim H^{1,0}(M_1)\geq \dim H^{1,0}(M_2)$.
Unless $\dim H^{1,0}(M_1)= \dim H^{1,0}(M_2)$, this would
imply that $\dim H^{1,0}(M_1) > \frac{b_1(M)-1}{2}$,
which is impossible because the map \eqref{_H^1_holom_map_Equation_}
is injective.

Notice that the equality $\dim H^{1,0}(M) = \frac{b_1(M)-1}{2}$
is valid for non-K\"ahler complex surfaces as well. 

The decomposition \eqref{_H^1_decompo_Equation_}
is the cornerstone for the description of the set of Lee classes
on an LCK manifold with potential. Consider the linear map
$\mu:\; H^1(M, \R)\arrow \R$ vanishing
on the codimension 1 subspace 
$H^{1,0}(M) \oplus \overline{H^{1,0}(M)}\subset H^1(M, \R)$
and positive on the Lee form. We prove that
$\xi \in H^1(M, \R)$ is a Lee class
if and only if $\mu(\xi) >0$
(\ref{_Lee_cone_on_LCK-pot_Theorem_}).

\hfill

\noindent{\bf Conventions:}  In the sequel, $(M,I)$ is a
connected complex manifold of complex dimension $n\geq
2$. For an Hermitian metric $g$, we shall denote with
$\omega(x,y):=g(Ix,y)$ the fundamental 2-form. We extend the
action of the complex structure to $k$-forms by
$(I\eta)(x_1,\ldots,x_k)=(-1)^k\eta(Ix_1,\ldots,Ix_k)$,
and we denote $I\eta$ by $\eta^c$. The complex
differential $d^c$ is defined as $d^c=I^{-1}dI$. We let
$d^*:\Lambda^kM\arrow \Lambda^{k-1}M$ be the metric adjoint
of the exterior derivative.

\hfill

\section{LCK manifolds}

We gather here the necessary background in LCK
geometry. For details, please see \cite{do} and
\cite{ov_jgp_16,ov_lckpot,ov_jgp_09}.

\subsection{Definitions. Examples}

\hfill

\definition\label{_LCK_def_via_formula_Definition_} 
$(M,I)$ is of {\bf locally conformally K\"ahler (LCK) type} if it admits an Hermitian metric $g$ whose fundamental form satisfies the equation
\begin{equation}\label{_def_LCK_equation_}
	d\omega=\theta \wedge \omega
\end{equation}
for a closed 1-form $\theta$ called {\bf the Lee form}. Then $(M,I,\omega,\theta)$ is called an {\bf LCK manifold}.

\hfill

\remark
\begin{enumerate}
	\item The LCK condition is conformally invariant: if $g$ is LCK with Lee form $\theta$, then $e^fg$ is LCK with Lee form $\theta+df$, hence to each conformal class of LCK metrics there corresponds a Lee class in $H^1(M,\R)$.
	\item In dimension $n\neq 2$, the equation \eqref{_def_LCK_equation_} implies $d\theta=0$.
	\item Using \eqref{_def_LCK_equation_}, one can
          prove that the Lee form is determined in terms
          of $I$ and $\omega$ by
          $\theta=-I\left(\frac{1}{n-1}d^*\omega\right)$.

\item If $\theta$ is exact, the LCK manifold is called
  {\bf globally conformally K\"ahler (GCK)}. Usually, it
  is tacitly assumed that $\theta$ is not exact.
\end{enumerate}

\hfill

In the sequel we will mostly use following definition, 
equivalent to \ref{_LCK_def_via_formula_Definition_}.

\hfill

\definition
A complex manifold $(M,I)$ is LCK if and only if it admits a cover $(\tilde M, I)$ equipped with a K\"ahler metric $\tilde\omega$ with respect to which the deck group of the cover acts by holomorphic homotheties.

\hfill

\definition The {\bf homothety character} associated to a K\"ahler cover with deck group $\Gamma$ is $\chi:\Gamma:\to\R^{>0}$, $\chi(\gamma)=\frac{\gamma^*\tilde\omega}{\tilde\omega}$.
The rank of $\Im(\chi)$ is the {\bf LCK rank} of $(M,I,\omega)$.

\hfill

\remark Since $\Gamma$ is a quotient group of $\pi_1(M)$,
we can consider $\chi$ as a character on $\pi_1(M)$. Let then $L\arrow M$ be the local system associated to $\chi$. It is a real line bundle and $\theta$ can be viewed as a connection form in $L$ which is thus flat. The line bundle $L$ is also called {\bf the weight bundle} of the LCK manifold.

\hfill

\definition {\bf The minimal K\"ahler cover} of an LCK
manifold corresponds to a group  $\Gamma$ on which $\chi$
is injective ($\Gamma$ does not contain
$\tilde\omega$-isometries). This is the smallest
cover admitting a K\"ahler metric which is
conformal to the pullback of the LCK metric.

\hfill

\definition 
A differential form $\alpha\in\Lambda^*\tilde M$ is called {\bf automorphic} if $\gamma^*\alpha=\chi(\gamma)\alpha$ for all $\gamma\in\Gamma$.

\hfill

\remark\label{_weight_bundle_remark_} 
\begin{description}
	\item[(i)] Automorphic forms on $\tilde M$ can be identified with $L$-valued forms on $M$. In particular, since $\pi^*\omega$ is $\Gamma$-invariant, $\omega$ can be viewed as a section of $\Lambda^{1,1}(M,L)$, and $\omega^k$ as a section of $\Lambda^{k,k}(M,L^{\otimes k})$ etc. 
	\item[(ii)] Let  $d_\theta:=d-\theta\wedge$. Then $d_\theta\omega=0$, hence $\omega$ is a closed $L$-valued form. 
	\item[(iii)] The complex $(\Lambda^*, d_\theta)$ is elliptic, since $d_\theta$ has the same symbol as $d$, and its cohomology $H^*_\theta(M)$ can be identified with the cohomology $H^*(M,L)$ of the local system $L$; it is called {\bf Morse-Novikov cohomology}.
\end{description}

\example
The following manifolds admit an LCK structure:
almost all known non-K\"ahler compact complex surfaces (see e.g. \cite{_ovv:surf_,va_tr, bel, go, _Brunella:Kato_}; Hopf manifolds: $(\C^n\setminus 0)/\langle A\rangle$, $A\in\mathrm{GL}(n,\C)$
with eigenvalues of absolute value $> 1$ (see e.g. \cite{ov_pams})\footnote{$(\C^n\setminus 0)/\langle A\rangle$ is called {\em diagonal Hopf manifold} when $A$ is diagonalizable, and {\em non-diagonal Hopf manifold} when $A$ is not diagonalizable.}; some Oeljeklaus-Toma manifolds (\cite{ot}); Kato manifolds (\cite{iop}) and some ``toric Kato manifolds'' (\cite{iopr}).

\subsection{The dichotomy K\"ahler {\em versus} LCK}

The next result, proven by Vaisman, shows that on compact complex manifolds, LCK and GCK metrics cannot coexist. For consistency, we provide a proof, slightly different from the original one.

\hfill  

\theorem {\rm (\bf \cite{va_tr})}\label{vailcknotk}
Let $(M,\omega, \theta)$ be a compact LCK manifold, not globally conformally K\"ahler. 
Then $M$ does not admit a K\"ahler structure.

\hfill

\proof {\bf Step 1:} That $M$ is not globally conformally K\"ahler means that $\theta$ is not cohomologuous with zero, that is $\theta$ is not $d$-exact. 

Let $d\omega=\omega\wedge\theta$,
$\theta'=\theta+d\phi$. Then
$$d(e^\phi\omega)= e^\phi\omega\wedge \theta+
e^\phi\omega\wedge d\phi =e^\phi\omega\wedge \theta'.$$
This means that
we can replace the triple $(M,\omega, \theta)$ by
$(M,e^\phi \omega, \theta')$ for any 1-form $\theta'$
cohomologous to $\theta$. 

\hfill

{\bf Step 2:} Assume that $M$ admits a K\"ahler structure.
Then, by Hodge theory, $\theta$ is cohomologous to the sum of a holomorphic and an antiholomorphic form.
After a conformal transformation (which changes $\theta$ in $\theta+d\phi$) as in Step 1,  we may assume
that $\theta$ itself  is the sum of a holomorphic and an  antiholomorphic form.

\hfill

{\bf  Step 3:} Then $dd^c\theta= \1 d \bar\6\theta=0$,
giving $dd^c(\omega^{n-1})= \omega^{n-1}\wedge \theta\wedge I(\theta)$. 
Therefore 
$0=\int_M dd^c(\omega^{n-1})=\int \mathrm{Mass}(\theta\wedge I(\theta))$,
hence $\theta\wedge I(\theta)=0$, thus $\Vert\theta\Vert^2=0$ and the 
initial metric is globally conformally K\"ahler.\footnote{Recall that the   mass of a positive $(1,1)$-form $\eta$, denoted $\mathrm{Mass}(\eta)$,  is the volume form
	$\eta \wedge \omega^{n-1}$, \cite{_Demailly:Book_}. }
\endproof

\hfill

Using similar techniques, we can prove:

\hfill

\lemma\label{_theta_not_d^c_closed_Lemma_}
Let $(M, \theta, \omega)$ be a compact LCK manifold. Then the cohomology class
$[\theta]\in H^1(M, \R)$ cannot be represented by a form which is
$d^c$-closed.

\hfill

\proof
Indeed, for each representative of $[\theta]$, this form can be realized
as the Lee form for an LCK metric which is conformally equivalent to $\omega$.
Therefore, it would suffice to show that $d^c\theta\neq 0$ for any compact LCK
manifold $(M, \theta, \omega)$. If $d^c\theta= 0$, we would have
$dd^c(\omega^{n-1})= \omega^{n-1}\wedge \theta\wedge I(\theta)$,
giving, as above, $0=\int_M dd^c(\omega^{n-1})=\int \mathrm{Mass}(\theta\wedge I(\theta))$,
hence $\theta\wedge I(\theta)=0$, implying that $\theta=0$.
\endproof

\hfill

This can be used to prove an important step in our
decomposition theorem for $H^1(M)$, where $M$ is an LCK
manifold with potential. 

\hfill 

\lemma\label{_H^1_holo_LCK_Lemma_}
Let $(M, \theta, \omega)$ be a compact LCK manifold,
and $H^{1,0}(M)$ denote the space of closed holomorphic 1-forms on $M$.
Then the natural map 
\[ H^{1,0}(M)\oplus \overline{H^{1,0}(M)} \oplus \langle \theta\rangle\arrow H^1(M,\C)
\]
is injective,
where $\langle \theta\rangle$ is the subspace generated by  $\theta$.

\hfill

\proof
A closed holomorphic form $\alpha$ belongs to $\ker d \cap \ker d^c$.
Indeed, $\bar\6 \alpha=0$ together with $d\alpha=0$ implies $d^c\alpha=0$.
Therefore, if  $\alpha\in H^{1,0}(M)+ \overline{H^{1,0}(M)}$
is exact, one has $\alpha = d f$ and $dd^c f=0$, which is impossible by the 
maximum principle. However, if $\theta$ is cohomologous to a 
sum of holomorphic and antiholomorphic forms, this easily leads
to a contradiction with \ref{_theta_not_d^c_closed_Lemma_}.
Indeed, suppose that $\theta=\alpha+ df$, where
$d\alpha=d^c\alpha=0$. Making a conformal change, we obtain
another LCK structure which has Lee form equal to $\alpha$.
This is impossible, again by \ref{_theta_not_d^c_closed_Lemma_}.
\endproof

\hfill

\corollary\label{_inequa_holo_LCK_Corollary_}
Let $M$ be a compact LCK manifold, and $H^{1,0}(M)$ denote the space of closed holomorphic 1-forms on $M$.
Then $\dim H^{1,0}(M) \leq \frac{b_1(M)-1}{2}$.
\endproof

%%%%%%%%%%%%%%%%%%%%%%%%%%%%%%%%%%%%%%%%%%%%%%%%%%%%%%%%%%%%%%%%%%
\section{Vaisman manifolds}\label{_Vaisman_manifolds_}
%%%%%%%%%%%%%%%%%%%%%%%%%%%%%%%%%%%%%%%%%%%%%%%%%%%%%%%%%%%%%%%%%%

The best understood subclass of LCK manifolds is the one
with the Lee form which is parallel with respect to the Levi-Civita
connection. They are called {\bf Vaisman
  manifolds}.

\hfill

If $(M,g,I,\theta)$ is Vaisman, the Lee field
$\theta^\sharp$ is Killing and holomorphic; moreover, it commutes with $I\theta^\sharp$ (\cite{va_gd}). Denote by $\Sigma$ the holomorphic 1-dimensional foliation generated by $\theta^\sharp$ and $I\theta^\sharp$. It is called {\bf the canonical foliation} (the motivation is given in the next theorem).

\hfill

\theorem \label{_Subva_Vaisman_Theorem_}
Let $M$ be a compact Vaisman manifold, and 
$\Sigma\subset TM$ its canonical foliation. Then:
\begin{description}
	\item[(i)] $\Sigma$ is
	independent from the choice of the Vaisman metric (\cite{tsu}).
	\item[(ii)] $d^c\theta=\omega-\theta \wedge I\theta$ (\cite{va_gd}) and the exact (1,1)-form $\omega_0:= d^c\theta$ is semi-positive
 (\cite{_Verbitsky:Vanishing_LCHK_}). Therefore, $\Sigma=\ker\omega_0$,
and $\omega_0$ is transversally K\"ahler with respect to $\Sigma$.
%	\item[(iii)] For any complex compact subvariety
%	$Z\subset M$, $Z$ is tangent to $\Sigma$. 
%	\item[(iv)] For any compact complex subvariety
%	$Z\subset M$, the set of smooth points of $Z$ is
%	Vaisman (\cite{_Verbitsky:Vanishing_LCHK_}).
\end{description}

Since $\theta^\sharp$ is Killing and holomorphic, it generates a complex flow of $g$-isometries. These  lift to
holomorphic non-trivial homotheties of the K\"ahler metric on
the  universal cover $\tilde M$. This is in
fact an equivalent definition of Vaisman-type manifolds,
as the following  criterion shows:

\hfill

\theorem{\bf  (\cite {kor})}\label{kami_or} 
Let $(M,\omega, \theta)$ be a compact LCK manifold equipped with a 
holomorphic and conformal $\C$-action $\rho$ without fixed points,
which lifts to non-isometric homotheties on 
the K\"ahler cover $\tilde M$. Then $(M,\omega, \theta)$
	is conformally equivalent to a Vaisman manifold.
	
	\hfill

\example\label{_Vaisman_Examples_}	
A non-exhaustive list of examples of Vaisman manifolds comprises:
\begin{description}
	\item[(i)] Diagonal Hopf manifolds $(\C^n\backslash 0)/\langle A\rangle$ where $A$ is semi-simple and with eigenvalues $\alpha_i$ of absolute value $>1$, \cite{go, ov_pams}. 
	\item[(ii)] Elliptic complex surfaces (see \cite{bel} for the complete classification of Vaisman compact surfaces; see also \cite{_ovv:surf_}). 
	\item[(iii)] All compact submanifolds of a Vaisman manifold (\cite{_Verbitsky:Vanishing_LCHK_}).
\end{description}

\remark The class of Vaisman manifolds is strict:  neither
the LCK Inoue surfaces, nor the non-diagonal
Hopf 
manifolds can bear Vaisman metrics (\cite{bel}, \cite{ov_pams}). 

\hfill

Recall that a form $\eta$ on a foliated manifold $(M, \Sigma)$
is called {\bf basic} if it can be locally obtained as a pullback
$\pi^*\eta_0$ from the leaf space of $\Sigma$, which is defined
in a sufficiently small neighbourhood of every point $x\in M$.

\hfill

The following claim is well known
(and can be used as a definition of basic forms).

\hfill

\claim\label{_basic_vanish_Claim_}
A form $\eta$ on $M$ is basic with respect to
$\Sigma\subset TM$ if and only if for any vector field
$X \in \Sigma$, one has $i_X(\eta) = \Lie_X(\eta)=0$,
where $i_X$ denotes the contraction with $X$.
\endproof

\hfill

\corollary\label{_closed_basic_Corollary_}
A closed form $\eta$ on $M$ is basic with respect to
$\Sigma\subset TM$ if and only if $i_X(\eta) = 0$.

\hfill

\proof
Follows from the Cartan formula $\Lie_X(\eta) = i_X(d\eta) + d(i_X(\eta))$.
\endproof

\hfill

Further on, we need the following observation.

\hfill

%%%%%%%%%%%%%%%%%%%%%%%%%%%%%%%%%%%%%%%%%%%%%%%%%%%%%%%%%%%%
\proposition\label{_holomo_on_Vaisman_basic_Proposition_}
Let $M$ be a compact Vaisman manifold, and $\eta$ a 
closed holomorphic 1-form on $M$. Then $\eta$ is basic
with respect to the canonical foliation $\Sigma$ on $M$.

\hfill

\proof
%Let $\theta^\sharp$ be the Lee field on $M$.
%The contraction $\langle \theta^\sharp, \eta\rangle$
%is a holomorphic function, because $\theta^\sharp$
%and $\eta$ are holomorphic. Since $M$ is compact,
%$\langle \theta^\sharp, \eta\rangle$ is constant. 
%By \ref{_closed_basic_Corollary_}, to prove that
%$\eta$ is basic, it would suffice to show that
%$\langle \theta^\sharp, \eta\rangle=0$.
Let $n=\dim_\C M$,
and $\omega_0\in \Lambda^{1,1}(M)$ 
the transversal K\"ahler form defined above.
Since $\eta$ is closed and $\omega_0$ is exact,
one has $\int_M \omega_0^{n-1}\wedge \eta\wedge\bar \eta=0$.
However, $-\1 \eta\wedge\bar \eta$
is a semi-positive form, and $\omega_0$
is strictly positive in the directions transversal
to $\Sigma$. This implies that
$-\1 \omega_0^{n-1}\wedge \eta\wedge\bar \eta$
is a positive volume form in every point
$x\in M$ such that $\eta\restrict {T_x M}$ does not vanish
on $\Sigma\restrict{T_x M}$. Since 
$\int_M \omega_0^{n-1}\wedge \eta\wedge\bar \eta=0$,
it follows that $\eta \restrict\Sigma=0$ everywhere.
By \ref{_closed_basic_Corollary_}, $\eta$ is basic.
\endproof

%%%%%%%%%%%%%%%%%%%%%%%%%%%%%%%%%%%%%%%%%%%%%%%%%%%%%%%%%%%%%%%%%%
\section{LCK manifolds with potential}\label{_LCK_w_potential_}
%%%%%%%%%%%%%%%%%%%%%%%%%%%%%%%%%%%%%%%%%%%%%%%%%%%%%%%%%%%%%%%%%%

We now introduce the main object of study of this paper.

\hfill

\definition An LCK manifold has  {\bf LCK potential} if it
admits a K\"ahler covering on which the K\"ahler metric
has a global and positive  potential function $\psi$ such that
the deck group multiplies $\psi$ by a constant. 
In this case, $M$ is called {\bf LCK manifold with
	potential}.

\hfill

\example All Vaisman manifolds are LCK manifolds with
potential. Indeed, if $\pi:\tilde M\arrow M$ is the universal cover and $\theta$ is the Lee form on $M$, then $\Vert\pi^*\theta\Vert$ is an automorphic global K\"ahler potential for $\tilde\omega$. Also, the structure is inherited by all complex submanifolds of an LCK manifold with potential. 
Among the non-Vaisman examples, we mention the
non-diagonal Hopf manifolds, \cite{ov_jgp_16}. 

\hfill

\remark
LCK Inoue surfaces, blow-ups of LCK manifolds 
and OT-manifolds cannot be LCK manifolds with potential, 
\cite{oti2,_Vuletescu:blowups_}.

\hfill

A wealth of examples is provided by the following fundamental result:

\hfill

\theorem {\bf (\cite{ov_lckpot})}
Let $(M,I,\omega,\theta)$ be a compact LCK manifold with
potential. Then any small deformation $(M,I_t)$, $t\in\C$,
$|t|<\epsilon$, admit an LCK metric with potential.
In particular, non-diagonal Hopf manifolds are LCK with potential.

\hfill

\remark \label{_d_theta_c_equation_}
By \cite{ov_imrn_10} (see also \cite{_Istrati:LCK-pot_}), 
an LCK manifold with potential admits a conformal gauge
such that
\begin{equation}\label{dctheta}
	d\theta^c=\omega - \theta\wedge\theta^c, \quad \text{where}\quad \theta^c(X)=-\theta(IX).
\end{equation}
We shall tacitly assume that the LCK metric is chosen in
such a way that \eqref{dctheta} holds.

\hfill

We are especially interested in the LCK manifolds with
potential of LCK rank 1, that is, LCK manifolds with
potential admitting a K\"ahler $\Z$-covering. 
We showed in \cite{ov_lckpot} that this is equivalent 
to the LCK potential being a proper function on
the minimal K\"ahler cover. The minimal cover of a
compact LCK manifold with proper potential is very
nice from the complex and algebraic viewpoint:

\hfill
 
\theorem {\bf (\cite{ov_lckpot,ov_pams})}\label{potcon}
Let $M$ be a compact LCK manifold with proper potential,
and $\tilde M$ its K\"ahler $\Z$-covering.
If  $\dim_\C M\geq 3$, then the metric completion $\tilde M_c$\index{completion!metric}
admits a structure of a complex variety, 
compatible with the complex structure on
$\tilde M \subset\tilde M_c$, and 
the complement $\tilde M_c\setminus  \tilde M$
is just one point. Moreover, $\tilde M_c$ is an affine
algebraic variety obtained as an affine cone
over a projective orbifold.

\hfill

\remark Notice that $\tilde M_c$ is indeed the {\bf Stein completion} of $\tilde M$ in the sense of \cite{andreotti_siu}. In the proof of our theorem, we used the filling theorem 
by Rossi and Andreotti-Siu (\cite{rossi,andreotti_siu}) which imposes the restriction $\dim_\C>2$. 

\hfill

By appearance, assuming that the potential is a proper
function is a restrictive condition. However, this is
not entirely true: as long as one
is interested in the complex geometry (and not in the
Riemannian one), one can always assume the LCK potential is proper.

\hfill 

\theorem {\em \bf(\cite{ov_jgp_09, ov_jgp_16})} 
\label{defor_improper_to_proper}
Let $(M, \omega, \theta, \phi)$ be a compact LCK manifold
with improper LCK potential. Then $(\omega, \theta, \phi)$
can be approximated in the ${C}^\infty$-topology by an LCK
structure with proper LCK potential on the same complex manifold.

\hfill

The following is one of the most important features of compact LCK manifolds with proper potential, of dimension greater than 3.

\hfill

\theorem\label{_Embedding_LCK_pot_in_Hopf_}
Any compact LCK manifold $M^n$, $n\geq 3$, admits a
holomorphic embedding into a Hopf manifold
$(\C^N\backslash 0)/\langle A\rangle$.
A manifold $M^n$ is Vaisman if
and only if it admits an embedding to
$(\C^N\backslash 0)/\langle A\rangle$, 
with the matrix $A$ diagonalizable.
 
\hfill

\proof
In \cite[Theorem 3.4]{ov_lckpot}
it is shown that an LCK manifold with potential
is embeddable into a Hopf manifold, and in
 \cite[Theorem 3.6]{ov_lckpot} it is shown
that a Vaisman manifold is embeddable to
a diagonal Hopf manifold. Conversely, in 
\cite[Section 2.5]{ov_pams} we show that
a diagonal Hopf manifold is Vaisman,
and in \cite{_Verbitsky:Vanishing_LCHK_} it is shown
that a positive-dimensional compact 
submanifold of a Vaisman manifold is Vaisman.
\endproof

\hfill

One of the most useful properties of compact LCK manifolds with potential in dimension greater than 3 is that their complex structure can be deformed to a complex structure that supports Vaisman metrics. 

\hfill

\theorem  {\bf (\cite{ov_imrn_10})} \label{def_lckpot2Vai}
Let $(M,\omega,\theta)$, 
$\dim_\C M \geq 3$, be a compact LCK manifold
with proper potential. Then there exists a complex
analytic deformation of $M$ which
admits a Vaisman metric.

\hfill

\remark A refinement of this result will be given in \ref{_Vaisman_limit_of_LCK_pot_Theorem_}.

\section{Algebraic cones and LCK manifolds with potential}

\subsection{Jordan-Chevalley decomposition}

%We recall the following definitions:
%
%\hfill
%
%\definition 
%\begin{enumerate}
%	\item The {\bf category of affine varieties} is
%	the category of finitely generated rings without
%	nilpotent elements  with arrows inverted.
%	\item An {\bf algebraic group} is a group object in the
%	category of affine varieties.
%	%	\item A {\bf pro-algebraic group} is an inverse limit
%	%of algebraic groups. 
%\end{enumerate}

Further on, all algebraic groups are
considered over $\C$. For the definition and more reference
on algebraic groups, please see \cite{hum}.

\hfill

\definition 
An element of an algebraic group $G$ is called
{\bf semisimple} if its image is semisimple
for some exact
algebraic representation of $G$, and is called 
{\bf unipotent} if its image is unipotent 
(that is, exponential of a nilpotent) 
for some exact algebraic representation of $G$.

\hfill

\remark 
For any algebraic representation of an algebraic group
$G$, the image of any semisimple element is a 
semisimple operator, and the image of any unipotent
element is a unipotent operator (\cite[\S 15.3]{hum}).

\hfill

\theorem  {\bf (Jordan-Chevalley  decomposition), 
	\cite[\S 15.3]{hum}} \label{jcdec}\\
Let $G$ be an algebraic group, and $A\in G$.
{Then there exists a unique decomposition $A= S U$ 
	of $A$ in a product of commuting elements $S$ and $U$,
	where $U$ is unipotent and $S$ semisimple.}

\hfill

%\remark 
%Since this decomposition is unique, it is functorial.
%{Therefore, it is also true for all pro-algebraic groups.}	

\subsection{Algebraic cones}

To better describe the universal cover of a compact LCK
manifold with potential we need to introduce the closed
and open algebraic cones.

\hfill

\definition 
A {\bf closed algebraic cone}  is an
affine variety $\cac$ admitting a $\C^*$-action $\rho$
with a unique fixed point $x_0$, called {\bf the origin},
and satisfying the following: 
\begin{description}
\item[(i)] $\cac$ is smooth outside of $x_0$,
\item[(ii)] $\rho$ acts on the Zariski tangent
space $T_{x_0}\cac$ with all eigenvalues
$|\alpha_i|<1$.
\end{description}

An {\bf open algebraic cone} is a closed algebraic
cone without the origin.

\hfill

%\remark
%It is not hard to see that an open algebraic
%cone is the total space of a principal
%$\C^*$-bundle $L_0$ over a projective orbifold,
%and the line bundle associated with $L_0$ is ample.

For the sake of completeness, we give a new and self-contained proof of the 
following basic result.

\hfill

\theorem {\bf (\cite[Theorem 2.8]{ov_pams})}
\label{_cone_cover_for_LCK_pot_Theorem_}
Let $M=\tilde M/\langle A\rangle$ 
be an LCK manifold with potential, with LCK rank 1, 
and  $\tilde M$  its K\"ahler $\Z$-covering. {Then $\tilde M$ is an open
	algebraic cone.}

\hfill

\pstep  Let $\tilde M_c$ be the Stein
completion of $\tilde M$ equipped with an $A$-equivariant embedding to
$\C^N$, where $A$ acts as a linear operator with all
eigenvalues $|\alpha_i|< 1$.
% (\ref{_gamma_on_T^*_contraction_Corollary_}). 

Let $\calo_{\C^N,0}$ denote the ring of germs of holomorphic functions in zero.
Call a function $f\in \calo_{\C^N,0}$  
{\bf $A$-finite} if the space $\langle f, A^*f, {A^2}^* f, ...\rangle$
is finite-dimen\-sio\-nal.
A polynomial function is clearly
$A$-finite. The converse is also true, because the Taylor 
decomposition of an $A$-finite function $f$ can  only have 
finitely many components, otherwise the eigenspace
decomposition of $f$ is infinite. 

\hfill

{\bf Step 2:}
We want to produce an explicit fundamental domain $U_0$
for the action of 
$\Z\cong \langle A\rangle =
\{ ..., A^{-n}, A^{-n+1}, ... A^{-1}, \Id_{\C^n}, A, A^2, ... \}$
on $\C^N$, in such a way that $U_0 = V \backslash A(V)$,
where $V$ is Stein. Let $B\subset \C^n$ be the unit ball. When the operator norm $\|A\|$ of $A$ 
is less than 1, one has $A(B) \Subset B$, and 
$B \backslash A(B)$ is the fundamental domain which we can use.
This would hold, for example, when $A$ is diagonalizable.
On the other hand, the operator norm of a contraction
can be bigger than 1. Consider for example 
the matrix
$A=\begin{pmatrix} \frac 1 2 & 1000 \\ 0 & \frac 1 2\end{pmatrix}$;
its norm is at least 1000. Therefore, one should take more
care when choosing the fundamental domain. Recall that any
matrix over $\C$ admits a Jordan decomposition, and
every Jordan cell 
{\begin{equation}\label{_Jordan_cell_Equation_}
\small 
	\begin{pmatrix} \alpha & 1 & 0 & \ldots & 0\\
		0 & \alpha & 1 & \ldots & 0\\
		\vdots &\vdots &\vdots & \cdots & \vdots \\
		0&0&0 & \ldots &1\\
		0&0&0 & \ldots &\alpha
	\end{pmatrix}
\end{equation}
}
is conjugate to 
{\begin{equation}\label{_Jordan_cell_epsilon_Equation_}
\small 
	\begin{pmatrix} \alpha & \epsilon & 0 & \ldots & 0\\
		0 & \alpha & \epsilon & \ldots & 0\\
		\vdots &\vdots &\vdots & \cdots & \vdots \\
		0&0&0 & \ldots &\epsilon\\
		0&0&0 & \ldots &\alpha
	\end{pmatrix}
\end{equation}}
(see {\em e. g.} \ref{_semisimple_operator_approx_Proposition_} below).
Writing $A$ in a Jordan basis and replacing each cell 
\eqref{_Jordan_cell_Equation_}
with \eqref{_Jordan_cell_epsilon_Equation_}, for $\epsilon$
sufficiently small, we obtain a contraction with an operator
norm $<1$ conjugate to $A$; then $A(B)\Subset B$
is a fundamental domain for the action of $\langle A\rangle$.

\hfill

{\bf Step 3:}
Let $U_0$ be a fundamental domain for $A$ acting
on $\C^N$. As indicated in Step 2,
every linear contraction is conjugate to an operator
$A$ with operator norm $<1$. 
The fundamental domain $U_0$ for $A$ with operator norm
$<1$ can be obtained by taking an open ball
$B\subset \C^n$ and removing $A(B)$ from $B$. Denote by
$U_n$ a copy of this domain obtained as $U_n:= A^{-n}(U_0)$,
and let $V_n:= \{0\}\cup \bigcup_{i>-\infty}^n U_i$.
Since  $V_n= A^{-n}(B)$, it is a Stein domain in $\C^N$.

Let $A:\; H^0_b(\calo_{V_n}) \arrow H^0_b(\calo_{V_n})$
be the operator on the ring of bounded holomorphic functions
induced by the action of $A$ on $\C^N$.
Clearly, this map is compatible with the map
$A:\; H^0(\calo_{\tilde M_c}) \arrow H^0(\calo_{\tilde M_c})$
constructed above; this is what allows us to denote
them by the same letter. 

\hfill

{\bf Step 4:}
We are going to prove that
the operator
\[ A:\; H^0_b(\calo_{V_n}) \arrow H^0_b(\calo_{V_n})
\]
is compact with respect to the topology defined by
the $\sup$-norm.\footnote{Recall that for a complex
  manifold $X$, the sup-topology on $H^0_b(\calo_X)$ is
  the topology given by the sup-norm, namely: $|f|_{\sup}
  := \sup_X |f|$.}

Let $\gamma:\; X \arrow X$ be a map
taking a complex manifold $X$ to its
precompact subset. We prove that in this case
the map $\gamma^*:\; H^0_b(\calo_{X}) \arrow H^0_b(\calo_{X})$
is compact in the $\sup$-topology.

For any $f\in H^0(\calo_X)$ we have
\[|\gamma^* f|_{\sup}= \sup_{x\in \overline{\gamma(X)}}
|f(x)|.
\]
This implies that $\gamma^*(f)$ is bounded.
Therefore, {for any sequence $\{f_i\in H^0(\calo_X)\}$ converging in the 
	$ C^0$-topology, the sequence $\{\gamma^* f_i\}$ converges
	in the $\sup$-topology.}
The set $B_C:=\{v\in V\ \ |
\ \  |v|_{\sup} \leq C\}$ is a normal family and hence, by 
Montel's theorem, it is precompact in the 
$ C^0$-topology (\cite[Chapter I, Theorem 3.12]{_Demailly:Book_}).
Then $\gamma^* B_C$ is precompact in the 
$\sup$-topology.
This proves that the operator $\gamma^*:\; V\arrow V$ is compact.

\hfill

{\bf Step 5:}
Let $\goth{I}(V_n)$ be the ideal of $\tilde M_c\cap V_n$
in $H^0_b(\calo_{V_n})$. Recall that $H^0_b(\calo_X)$
is a Banach algebra, by Montel's theorem (\cite[Chapter
IX, Proposition 4.7]{_Demailly:Book_}).
By the Riesz-Schauder's theorem (\cite[Section
  5.2]{friedman}), a compact endomorphism of 
a Banach space admits a Jordan decomposition.
Then $A$-finite vectors are finite linear combinations
of the vectors from the Jordan cells. This implies that
the set of $A$-finite 
functions in $\goth{I}(V_n)$
is dense in $\goth{I}(V_n)$, with respect to the $\sup$-topology.
On the other hand all $A$-finite functions can be holomorphically
extended to $\C^N$ by automorphicity. 

The base of $C^0$ (that is, compact-open) topology on $H^0(\calo_{\C^N})$ 
is formed by translations of open sets consisting of all
functions $f\in H^0(\calo_{\C^N})$ 
which satisfy $|f|< C< \infty$ on a given compact, for some
positive $C\in \R$. Therefore,
it is the weakest topology such that
its restriction to $H^0_b(\calo_{V_n})$
with $\sup$-topology is continuous.
This implies that any set of functions ${\goth S}\subset H^0(\calo_{\C^N})$
which is bounded on compacts and dense in $H^0_b(\calo_{V_n})$,
for all $n$, is dense in $H^0(\calo_{\C^N})$.

Since the space of $A$-finite functions
is dense in $\goth{I}(V_n)$,
the space $\goth{I}^A$  of $A$-finite functions in 
$\goth{I}$ is dense in $\goth{I}$ with respect to the
$C^0$-topology. In particular,
the set of common zeros of $\goth{I}^A$
coincides with $\tilde M_c\subset \C^N$.

\hfill

{\bf Step 6:}
The $A$-finite functions are polynomials, as shown in Step 1.
By Hilbert's basis theorem,
any ideal in the ring of polynomials is finitely generated.
Therefore, the ideal $\goth{I}^A$ is finitely generated
over polynomials. Let $f_1, ..., f_n$ be the
set of its generators. By Step 2, the set of
common zeros of $\goth{I}^A$ is 
$\tilde M_c\subset \C^N$; therefore,
$\tilde M_c\subset \C^N$ is given by  polynomial
equations $f_1=0, f_2=0, ...,  f_n=0$.

\hfill

{\bf Step 7:} It remains to show that
$\tilde M_c$ admits a holomorphic $\C^*$-action 
containing a contraction.
Let $G$ be the Zariski closure of
$\langle A \rangle$ in $\GL(\C^N)$. This is a commutative
algebraic group, acting on the variety $\tilde M_c\subset \C^N$.
Let $A=SU$ be the Jordan-Chevalley decomposition for
$A$, with $S, U\in G$. Since $G$ preserves $\tilde M_c$, 
theendomorphisms $S$ and $U\in\End\C^n)$ also act on $\tilde M_c$.
Since the eigenvalues of $S$ are the same
as the eigenvalues of $A$, it is a contraction.
Let $G_S\subset G$, $G_S= e^{\C \log S}$ be a one-parametric subgroup
containing $S$. We prove that $G_S$ can be approximated by
subgroups of $G$ isomorphic to $\C^*$; then these subgroups
also contain a contraction, and we are done.

Consider the map taking any $A_1\in G$ to its unipotent component
$U_1$. Since $G$ is commutative, this map is a group
homomorphism. Therefore, its kernel $G_s$ (that is,
the set of all semisimple elements in $G$) is an
algebraic subgroup of $G$. A semisimple commutative
algebraic subgroup of $\GL(\C^N)$ is always isomorphic
to $(\C^*)^k$ (\cite[Proposition 1.5]{_Borel_Tits:Groupes_Reductifs_}).
The one-parametric subgroups $\C^*\subset (\C^*)^k$ 
are dense in $(\C^*)^k$ because one-parametric
complex subgroups $\C^*\subset (\C^*)^k$ can be obtained as 
complexification of subgroups $S^1\subset U(1)^k \subset (\C^*)^k$,
and those are dense in $U(1)^k$. Therefore,
the contraction $S\in G_s=(\C^*)^k$ can be approximated
by an element of $\C^*$ acting on $\tilde M_c$.
\endproof

\hfill

\section{Hodge decomposition for $H^1(M)$ on LCK manifolds with potential}
%%%%%%%%%%%%%%%%%%%%%%%%%%%%%%%%%%%%%%%%%%%%%%%%%%%%%%%%%%%%%%%%%%%%%%%%

Any harmonic $r$-form, $r\leq n-1$,  on a compact $n$-dimensional Vaisman manifold $(M,\omega,\theta)$ can be uniquely written as a sum $\alpha+\theta\wedge\beta$ where $\alpha$ and $\beta$ are basic harmonic forms (see \cite{va_gd} or \cite{_ov_super_sas_} for a different proof. In particular, the space of harmonic 1-forms on a compact Vaisman manifolds is identified with $\ker d \cap \ker d^c\oplus \langle \theta\rangle$.

For LCK manifolds with potential, such a decomposition is no longer available. Instead we can prove:

\hfill

%%%%%%%%%%%%%%%%%%%%%%%%%%%%%%%%%%%%%%%%%%%%%%%%%%%%%%%%%%%%
\theorem\label{_LCK_pot_Hodge_decompo_Theorem_}
Let $(M, \theta, \omega)$ be a compact LCK manifold with potential,
and $H^{1,0}(M)$ denote the space of closed holomorphic 1-forms
on $M$. Using \ref{_H^1_holo_LCK_Lemma_}, we consider $H^{1,0}(M)\oplus
\overline{H^{1,0}(M)}\oplus \langle \theta \rangle$ as a subspace in $H^1(M, \C)$.
Then $H^1(M,\C) = H^{1,0}(M)\oplus \overline{H^{1,0}(M)}
\oplus \langle \theta \rangle$.

\hfill

\proof
%For Vaisman manifolds, this result is already proven 
%(\ref{_harmo_deco_1-form_Proposition_}). Indeed, by \ref{_harmo_deco_1-form_Proposition_}
%the space of harmonic forms in $\Lambda^1(M,\C)$ is 
%identified with $\ker d \cap \ker d^c\oplus \langle \theta\rangle$.
%The injectivity of the map 
%\[ H^{1,0}(M)\oplus \overline{H^{1,0}(M)} \oplus \langle \theta\rangle
%\arrow H^1(M,\C)
%\]
%follows from \ref{_H^1_holo_LCK_Lemma_}.
To prove that the map
\[ H^{1,0}(M)\oplus \overline{H^{1,0}(M)} \oplus \langle \theta\rangle
\arrow H^1(M,\C)
\]
is surjective, it would suffice to show that
$\dim_\C H^{1,0}(M)= \frac{b_1(M)-1}{2}$. We prove it by deforming
$M$ to a Vaisman manifold $M_0$ and showing that 
$\dim H^{1,0}(M)= \dim H^{1,0}(M_0)$.

We first deform the LCK metric on $M$ to an LCK metric
of LCK rank 1 (\ref{defor_improper_to_proper}). This operation does not affect the complex structure
on $M$, hence $\dim H^{1,0}(M)$ does not change, and it will suffice to 
prove that $\dim_\C H^{1,0}(M)= \frac{b_1(M)-1}{2}$ when 
$M$ is an LCK manifold with proper potential.

Let $\tilde M$ be the open algebraic cone associated with $M$ as in 
\ref{_cone_cover_for_LCK_pot_Theorem_}, and $A:\; \tilde M \arrow \tilde M$
the generator of the deck group. Applying the Jordan-Chevalley 
decomposition $A=SU$ as in \ref{def_lckpot2Vai}, we can deform
$\tilde M/\langle A\rangle$ to the Vaisman manifold
$M_0:=\tilde M/\langle S\rangle$. To prove that 
$\dim_\C H^{1,0}(M) = \frac{b_1(M)-1}{2}$ it would suffice to show that
all holomorphic, $S$-invariant 1-forms on $\tilde M$ are
also $U$-invariant. 

Consider $U$ as an automorphism of $M_0$. This automorphism
is homotopy equivalent to the identity because
$U= e^{N}$, where $N$ commutes with $S$.
Since $U$ is an unipotent element of the group
of automorphisms of the algebraic cone $\tilde M$,
the action of $U_t:= e^{tN}$ preserves $\tilde M$ and 
commutes with $S$, hence it is well defined on $M_0$.
This gives a homotopy of $U=U_1$ to $\Id=U_0$.

Since $U$ is homotopy equivalent to the identity,
it acts trivially on $H^1(M_0)$, hence
all $S$-invariant holomorphic forms
on $\tilde M$ are also $SU$-invariant.
This implies that $\dim_\C H^{1,0}(M) \geq \dim_\C H^{1,0}(M_0)  =\frac{b_1(M_0)-1}{2}$.
The inequality in this expression is in fact an equality by
\ref{_inequa_holo_LCK_Corollary_}. 
We thus proved \ref{_LCK_pot_Hodge_decompo_Theorem_}.
\endproof

\section{Approximating LCK with potential structures  by Vaisman structures}

We start with a linear algebra result which
will be used in the proof of the main theorem of this
section:

\hfill

\proposition\label{_semisimple_operator_approx_Proposition_}
Let $p\in \GL(n,\C)$ be a linear operator, and
$p=su$ its Jordan decomposition, with $s$ semisimple,
$u$ unipotent, and $su=us$. Then there exists
a sequence $r_i\in \GL(n,\C)$ of operators
commuting with $s$ and satisfying
$\lim_{i\to\infty} r_i p r_i^{-1}=s$.

\hfill

\proof
Since any operator is a sum of Jordan cells,
it would suffice to prove
\ref{_semisimple_operator_approx_Proposition_}
when $p$ is a single $k\times k$ Jordan cell,
{\[\small 
	p =\begin{pmatrix} \alpha & 1 & 0 & \ldots & 0\\
		0 & \alpha & 1 & \ldots & 0\\
		\vdots &\vdots &\vdots & \cdots & \vdots \\
		0&0&0 & \ldots &1\\
		0&0&0 & \ldots &\alpha
	\end{pmatrix}
	\]}
In this case, $s= \const\Id$, hence
it commutes with everything. Take
{\[\small 
	r_i =\begin{pmatrix} 1 & 0 & 0 & \ldots & 0\\
		0 & \epsilon_i & 0 & \ldots & 0\\
		\vdots &\vdots &\vdots & \cdots & \vdots \\
		0&0&0 & \ldots &0\\
		0&0&0 & \ldots &\epsilon_i^k
	\end{pmatrix}
	\]}
Then 
{\[\small 
	r_i p r_i^{-1} =\begin{pmatrix} \alpha & \epsilon_i & 0 & \ldots & 0\\
		0 & \alpha & \epsilon_i & \ldots & 0\\
		\vdots &\vdots &\vdots & \cdots & \vdots \\
		0&0&0 & \ldots &\epsilon_i\\
		0&0&0 & \ldots &\alpha
	\end{pmatrix}
	\]}
Taking a sequence $\epsilon_i$ converging to 0,
we obtain $\lim_{i\to\infty} r_i p r_i^{-1}=s$.
\endproof 

\hfill

The main  result of this section, \ref{_Vaisman_limit_of_LCK_pot_Theorem_},  
gives a more precise description of the approximation 
 in \ref{def_lckpot2Vai}. In order to state it, we need to recall the notion of Teichm\"uller space. Recall first that the {\bf $C^k$-topology}
on the space of sections of a bundle $B\arrow M$
is the topology of uniform convergence 
of $b$, $\nabla b$, $\nabla^2 b$, ...,
$\nabla^k b$ on compacts, for some
connection $\nabla$ on $B$. 
The {\bf $C^\infty$-topology} 
is the topology of uniform convergence
of {\em all} derivatives. In other words,
a set is open in the $C^\infty$-topology
if it is open in all $C^k$-topologies. Now we can give:

\hfill

\definition
Let $\Comp$ be the set of all integrable complex
structures on $M$, equipped with the $C^\infty$-topology,
and $\Diff_0$ the group of isotopies of $M$, that is, the
connected component of the group of diffeomorphisms of $M$.
{\bf The Teichm\"uller space} of complex structures
on $M$ is the quotient $\Teich:=\Comp/\Diff_0$
equipped with the quotient topology.

\hfill

%%%%%%%%%%%%%%%%%%%%%%%%%%%%%%%%%%%%%%%%%%%%%%%%%%%%%%%%%%%%
\theorem\label{_Vaisman_limit_of_LCK_pot_Theorem_}
Let $(M, J)$ be an LCK manifold with potential, $\dim_\C M \geq 3$.
Then there exists a Vaisman-type complex structure $(M, J_\infty)$
such that the point $[J_\infty]$ in the Teichm\"uller
space $\Teich(M)$ of complex structures on $M$ belongs to the 
closure of $[J]\in \Teich(M)$. In other words, there exists
a sequence of diffeomorphisms $\nu_i \in \Diff_0(M)$
such that $\lim_i \nu_i(J)=J_\infty$, where the limit is taken
with respect to the $C^\infty$-topology on 
the space $\Comp$ of complex structures.

\hfill

\proof
Without restricting the generality, we may
assume that $(M,J)$ has LCK rank 1, and its
$\Z$-cover is K\"ahler (\ref{defor_improper_to_proper}). 
Fix an embedding of $(M,J)$ into the Hopf manifold 
$H= \frac{\C^n\backslash 0}{\langle A \rangle}$ (\ref{_Embedding_LCK_pot_in_Hopf_}).

Let $A=us$ be the Jordan-Chevalley decomposition
for $A$. By \ref{_semisimple_operator_approx_Proposition_},
there exists a sequence $A_i= u_i s= r_i A r_i^{-1}$ 
of operators conjugated to $A$ such that $u_i$ converges
to 0. Denote by $(H, I_i)$ the Hopf manifold 
$(H, I_i):= \frac{\C^n\backslash 0}{\langle A_i \rangle}$.
Since $(H, I_i)$ are all naturally 
isomorphic to $H$, one obtains the embedding
$\phi_i:\; (M,J)\arrow (H, I_i)$.

Since the operators $A_i= u_i s$ converge 
to $s$, the sequence $I_i\in \Comp(H)$
converges to $I_\infty$, where
$(H, I_\infty):= \frac{\C^n\backslash 0}{\langle s \rangle}$.

 Denote by
$\gamma$ the generator of the monodromy acting on the Stein completion $\tilde M_c$ (\ref{potcon}), 
and $\phi:\; \tilde M_c\arrow \C^n$
the embedding making the following diagram
commutative
\[
\begin{CD}
	\tilde M_c@>{\phi}>> \C^n\\
	@V{\gamma}VV @V{A}VV\\
	\tilde M_c@>{\phi}>> \C^n.
\end{CD}
\]
Consider the map $\phi_i := r_i \phi r_i^{-1}$
as an embedding from $\tilde M_c$ to $\C^n$
making the following diagram commutative
\[
\begin{CD}
	\tilde M_c@>{\phi_i}>> \C^n\\
	@V{\gamma}VV @V{A_i}VV\\
	\tilde M_c@>{\phi_i}>> \C^n.
\end{CD}
\]
We use the same letter $\phi_i$ to denote
the embedding 
$(M, J) \arrow \frac{\C^n\backslash 0}{\langle A_i  \rangle}$ associated with $\phi_i$.
Since $\phi(\tilde M_c)$ is $s$-invariant, and $A_i$
converge to $s$, the
sequence $\phi_i\restrict{\tilde M_c}$ converges to 
$\phi\restrict {\tilde M_c}$, giving
an embedding 
$\frac{\tilde M}{\langle s\rangle} \arrow
(H,I_\infty)$. The  limit manifold 
$(M, J_\infty)=\frac{\tilde M}{\langle s\rangle}$ is of Vaisman type,
because it is embedded to a diagonal Hopf manifold
(\ref{_Embedding_LCK_pot_in_Hopf_}).

The maps $\phi_i$ do not converge to $\phi_\infty$ smoothly,
because the sequence $\{r_i^{-1}\}$ is not bounded.
However, the sequence $\{(\phi_i(\tilde M_c), A_i)\}$
$C^\infty$-converges to $(\phi_\infty(\tilde M_c), A)$ as a sequence
of pairs 
\[ \text{(algebraic subvariety $Z\subset \C^n$,  an
automorphism $\psi\in \Aut(Z)$);}
\] hence the
corresponding points in $\Teich$ also converge.
This is what we are going to show.

Let $S\subset \calo_{\C^n}$, $\dim S=m$, be a finite-dimensional space
generating the ideal of $\phi_1(\tilde M_c)$
(\ref{_cone_cover_for_LCK_pot_Theorem_}).
By Step 1 of the proof of 
\ref{_cone_cover_for_LCK_pot_Theorem_}, we may assume that all elements of $S$
are %homogeneous 
polynomials of degree less than $d$. 
Denote by $V\subset \calo_{C^n}$ the space of
polynomials of degree $\leq d$, and let $X\subset \Gr_m(V)$ be 
a subset of the Grassmannian of $m$-dimensional planes in $V$
consisting of all subspaces $W\subset V$ which generate
an ideal $J_W\subset \C[\C^n]$ in the polynomial ring
such that $\C[\C^n]/J_W$ is isomorphic to the algebraic cone 
$\tilde M_c$.  \footnote{Interpreting $X$ as a piece of the
relevant Hilbert scheme, we obtain that it is an
algebraic subvariety in $\Gr_m(V)$; we are not going to prove or use
this observation.}

The sequence
$\{\sigma_i :=\phi_i(\tilde M_c)\}$ corresponds to
points in $X$ converging to $\sigma_\infty:=\phi_\infty(\tilde M_c)$.
This gives the convergence of the
submanifolds $\phi_i(\tilde M)$ to $\phi_\infty(\tilde M)$
in $\C^n \backslash 0$.
Indeed, consider the ``universal
fibration'' over $X$, with the fiber over $W\in X$ 
being the algebraic cone associated with the ideal
$J_W\subset \calo_{\C^n}$ generated by $W$. 
The associated open cone fibration 
has smooth fibers. Indeed, any open algebraic cone 
is the total space of a $\C^*$-bundle over
a projective manifold.

To finish the proof, we need to prove that
the manifolds $(M,J_i)=\frac{\phi_i(\tilde M)}{\langle u_i s\rangle}$
smoothly converge to 
$(M, J_\infty)=  \frac{\phi(\tilde M)}{\langle s\rangle}$.
This would follow if we prove that 
the corresponding cones in $\C^n$ converge smoothly in each
annulus $B_R \backslash B_r$  around 0 (we need to restrict to the annulus,
because the cone itself is singular around zero, hence it makes
no sense to speak of $C^\infty$-convergence unless we remove a
neighbourhood of the origin).
Then $(M, J_i)$ and $(M, J_\infty)$ are quotients
of the respective cones by $A_i$ and $A$-actions
respectively, and $A_i$ converges to $A$ in $\GL(n, \C)$.

However, the cones $\phi_i(\tilde M_c)$ are smooth in each annulus, and they
converge to $\tilde M_c$ in $C^0$-topology (or in the Hausdorff metric)
by construction. For smooth families of compact
manifolds, the  $C^\infty$-convergence of their fibers is automatic.
To finish the proof, we replace the cone fibration 
over $X$
by the corresponding fibration of compact complex orbifolds,
which also converges to the central fiber. The fibers
of a locally trivial fibration of compact orbifolds
converge to the central fiber with all derivatives
by Ehresmann's theorem.

Let $P\arrow X$
be the fibration with 
projective fibers over $X$, obtained by
taking $\C^*$-quotients of the tautological open cone
fibration $U\arrow X$. 
The fibration $P\arrow X$ is locally trivial,
because it is smooth and all its fibers are isomorphic
projective orbifolds. The fibers of $U$
are total spaces of the $\C^*$-bundles associated
with $\calo(1)$ over fibers of $P$.
Then $U\arrow X$ is smoothly locally trivial.

To obtain the convergence of the corresponding
LCK manifolds, we notice that  $\lim_i A_i=s$,
hence $\lim_i \tilde M_i /\langle A_i \rangle = 
\tilde M /\langle s \rangle= (M, J_\infty)$.
\endproof

\hfill

\corollary
Let $(M,J)$ be a compact complex manifold, $\dim_\C M\geq 3$,
admitting an LCK structure with potential, and
$J_\infty$ the Vaisman-type complex structure
on $M$ obtained as in \ref{_Vaisman_limit_of_LCK_pot_Theorem_}.
Then any Vaisman-type Lee form on $(M, J_\infty)$
can be realized as the Lee form of an LCK
structure with potential on $(M, J)$.

\hfill

\proof Let $(M, I, \omega, \theta)$ be a Vaisman structure on $M$,
with $\omega= d^c \theta + \theta \wedge \theta^c$, and
$I_i$ a sequence of complex structures on $M$ converging
to $I$, such that all $(M,I_i)$ are isotopic to  $(M,J)$ 
as complex manifolds. Then the sequence
$\omega_i= I_i d I_i^{-1} \theta + \theta \wedge I_i\theta$
converges to 
\[ d^c \theta + \theta \wedge \theta^c=
I d I^{-1} \theta + \theta \wedge I\theta.
\]
Since positivity is an open condition, the (1,1)-form
$\omega_i$ is positive for $i$ sufficiently big.
Then $(M, I_i, \omega_i, \theta)$ is LCK with potential,
and $\theta$ its Lee form. However, $I_i$ is mapped
to $J$ by an isotopy which preserves the cohomology
class of $\theta$, hence $\theta$ is a Lee class on $(M, J)$.
\endproof

%%%%%%%%%%%%%%%%%%%%%%%%%%%%%%%%%%%%%%%%%%%%%%%%%%%%%%%%%%%%%%%%%%%%%%%%
\section{The set of Lee classes}
%%%%%%%%%%%%%%%%%%%%%%%%%%%%%%%%%%%%%%%%%%%%%%%%%%%%%%%%%%%%%%%%%%%%%%%%
%%%%%%%%%%%%%%%%%%%%%%%%%%%%%%%%%%%%%%%%%%%%%%%%%%%
\subsection{Opposite Lee forms on LCK manifolds with potential}
%%%%%%%%%%%%%%%%%%%%%%%%%%%%%%%%%%%%%%%%%%%%%%%%

As another preliminary result, we need
the following non-existence claim. For Vaisman manifolds,
it was obtained by K. Tsukada (\cite{tsuk}).

\hfill

%%%%%%%%%%%%%%%%%%%%%%%%%%%%%%%%%%%%%%%%%%%%%%%%
\proposition\label{_Lee_cannot_be_opposite_Proposition_}
Let $(M, \theta, \omega)$ and $(M, \theta_1, \omega_1)$
be two LCK structures on the same
compact complex manifold. Suppose that
$(M, \theta, \omega)$ is an LCK structure with potential.
Then $\theta+\theta_1$ 
cannot be cohomologous to 0.

\hfill

\pstep
If $[\theta]$ is the Lee class for an 
LCK structure with potential on $M$, then
$a[\theta]$ is also a Lee class for one, for any
$a>1$. To see this, consider 
the expression 
$\omega=d^c \theta + \theta \wedge\theta^c$ (\ref{_d_theta_c_equation_}) corresponding to the K\"ahler potential
$\phi$ on the K\"ahler cover $(\tilde M, \tilde \omega)\arrow (M,\theta,\omega)$,
with $\pi^*\theta = -d\log\phi$.
Then $\phi^a$ is also a K\"ahler potential
on $\tilde M$,
\[ 
dd^c \phi^{a} = \phi^{a-2} (a \cdot \phi dd^c \phi + a(a-1)
d\phi\wedge d^c\phi).
\]
Indeed, the first summand $a\phi^{a-1} dd^c \phi$ is Hermitian, because
$dd^c\phi$ is Hermitian, and the second summand
$a(a-1)d\phi\wedge d^c\phi$ is positive.
The function $\phi^a$ is automorphic,
hence it defines an LCK structure with
potential on $M$, and the corresponding
Lee form is  $- d \log (\phi^a)=a\theta$.

\hfill

{\bf Step 2:}
Let $\omega, \omega_1$ be LCK forms, and
$\theta, \theta_1$ the corresponding
Lee forms. Suppose that $k \theta + l \theta_1=0$.
Then 
\begin{equation*}%\label{_omega_k_wedge_omega_l_closed_equation_}
	 d(\omega^k \wedge \omega_1^l)= 
d\omega^k \wedge \omega_1^l + \omega^k \wedge
d\omega_1^l= k\theta \wedge \omega^k \wedge \omega_1^l
+ l\theta_1 \wedge \omega^k \wedge \omega_1^l=0.
\end{equation*}
This computation can be interpreted as follows.
Let $L$ be the weight bundle for $(M, \omega, \theta)$
and $L_1$ the weight bundle for $(M, \omega_1, \theta_1)$.
Recall (\ref{_weight_bundle_remark_} (ii)) that $\omega$, $\omega_1$  are viewed as closed $L$- and
$L_1$-valued forms. Then 
 $\omega^k$ is a closed $L^{\otimes k}$-valued form,
$\omega_1^l$ is a closed $L_1^{\otimes l}$-form,
and $\omega^k \wedge \omega_1^l$ is a closed
form with coefficients in the flat bundle
$L^{\otimes k}\otimes L_1^{\otimes l}$,
which is trivial.

Return now to the situation described in 
the assumptions of
\ref{_Lee_cannot_be_opposite_Proposition_}.
Let $n=\dim_\C M$.
Using Step 1, we replace the LCK structure
$(\omega, \theta)$ by another LCK structure with potential
in such a way that $\theta$ is replaced by
$(n-1)\theta$. Then $(n-1)\theta_1 = -\theta$,
and the volume form
$\omega \wedge \omega_1^{n-1}$ is closed.
However, $\omega$ is actually an exact $L$-valued
form, because $\omega= d_\theta (\theta^c)$, hence
$\omega \wedge \omega_1^{n-1}$ is an exact
$L \otimes L_1^{\otimes (n-1)}$-valued form.
However, $L \otimes L_1^{\otimes (n-1)}$ is a trivial
local system, which implies that $\omega \wedge \omega_1^{n-1}$ 
is exact.

We verify this with an explicit computation:
\begin{equation*}
	\begin{split}
		d(\theta^c \wedge \omega_1^{n-1})&=
		d\theta^c \wedge \omega_1^{n-1}-
		\theta^c \wedge d\omega_1^{n-1}\\
		&=(\omega - \theta\wedge \theta^c)\wedge \omega_1^{n-1}
		- (n-1)\theta_1 \wedge \theta^c \wedge \omega_1^{n-1}\\
		&=\omega \wedge \omega_1^{n-1} -
		(\theta\wedge \theta^c+(n-1)\theta_1 \wedge
		\theta^c)\wedge \omega_1^{n-1}=
		\omega \wedge \omega_1^{n-1}.
	\end{split}
\end{equation*}
We have shown that the positive volume form
$\omega \wedge \omega_1^{n-1}$ on $M$ is exact,
which is impossible.
\endproof

%%%%%%%%%%%%%%%%%%%%%%%%%%%%%%%%%%%%%%%%%%%%%%%%%%%%%%%%%%%%%%%%%%%%%%%%
\subsection{The set of Lee classes on Vaisman manifolds}
%%%%%%%%%%%%%%%%%%%%%%%%%%%%%%%%%%%%%%%%%%%%%%%%%%%%%%%%%%%%%%%%%%%%%%%%

To proceed, we need the following preliminary result,
which might be of separate interest.

\hfill

%%%%%%%%%%%%%%%%%%%%%%%%%%%%%%%%%%%%%%%%%%%%%%%%
\proposition\label{_Lee_form_on_Vaisman_is_Vaisman_Proposition_}
Let $(M,\theta, \omega)$ be an LCK structure
on a compact Vaisman manifold. Then $\theta$
is cohomologous to a Lee form of a Vaisman structure.

\hfill

\proof
Let $X$ be the Lee field of a Vaisman structure $(M, \omega^V, \theta^V)$
on $M$, and $G$ the closure of the group generated by 
exponents of $X$ and $I(X)$. Since $X$ and
$IX$ are Killing and commute, $G$ is a compact
commutative Lie group, hence it is isomorphic to a compact torus.
This group acts on $M$ by holomorphic isometries
with respect to the Vaisman metric.

Averaging $\theta$ with the $G$-action, we obtain
a $G$-invariant 1-form $\theta$, corresponding
to another LCK structure in the same conformal class.
Without restricting the generality, we may assume
from the beginning that the form $\theta$ is $G$-invariant. 

Now, the equation
$d\omega=\omega\wedge\theta$ is invariant under the
action of $G$, because $\theta$ is $G$-invariant;
in other words,
$d(g^*\omega)= g^*\omega\wedge \theta$, for all $g\in G$.
This implies that
$\omega$ averaged with $G$ gives a form
$\omega^G$ which satisfies
$d(\omega^G)= \omega^G\wedge \theta$.
We have constructed a $G$-invariant
LCK structure $(M, \omega^G, \theta)$. 

After lifting it to the K\"ahler
cover $\tilde M$ of $(M, \omega^G, \theta)$, the group $G$ becomes
non-compact. Indeed, if it remained compact,
it would act by isometries on the universal
cover $\tilde M_U$ of $M$ as well, hence 
the action of $G$ on $M$ is lifted to 
the action of $G$ on $\tilde M_U$.
This is impossible, however,
because the lift of $G$ to the
K\"ahler cover associated with $(M, \omega^V, \theta^V)$
acts by non-trivial homotheties, hence $G$
is lifted to an infinite cover $\tilde G\arrow G$
effectively acting on $\tilde M_U$. 
%(\ref{_Vaisman_Lee_action_contains_monodromy_Corollary_}). 

We obtain that $G$ 
acts by non-isometric homotheties on 
the K\"ahler cover associated with
$(M, \omega^G, \theta)$.
By 
\ref{kami_or}, $(M, \omega^G, \theta)$ is actually
Vaisman.
\endproof

\hfill

The following result was obtained by K. Tsukada
(\cite{tsuk}). We provide a new, simpler proof.

\hfill

%%%%%%%%%%%%%%%%%%%%%%%%%%%%%%%%%%%%%%%%%%%%%%%%%%%%%%%%%%%%
\theorem \label{_Lee_cone_on_Vaisman_Theorem_}
Let $M$ be a compact Vaisman manifold, and
$H^1(M)= H^{1,0}(M) \oplus \overline{H^{1,0}(M)} \oplus \langle \theta\rangle$
be the decomposition established in \ref{_LCK_pot_Hodge_decompo_Theorem_}.
Consider a 1-form $\mu\in H^1(M)^*$ vanishing on 
$H^{1,0}(M) \oplus \overline{H^{1,0}(M)} \subset H^1(M)$
and satisfying $\mu([\theta])>0$. Then 
a class $\alpha\in H^1(M,\R)$ is a 
Lee class for some LCK structure if and only if $\mu(x) >0$.

\hfill

\pstep
We start by proving that any $\alpha\in H^1(M,\R)$ satisfying
$\mu(x) >0$ can be realized as a Lee class. 

From \ref{_d_theta_c_equation_},
we have $\omega=d^c\theta+\theta\wedge I\theta$.
By \ref{_Subva_Vaisman_Theorem_} (ii),  the form
$\omega_0:=d^c\theta$ is semi-positive: 
it vanishes on the canonical foliation $\Sigma$ 
and is strictly positive in the transversal directions. 
Let $u\in H^{1,0}(M) \oplus \overline{H^{1,0}(M)}$.
Then
\[
d^c(\theta+ u) + (\theta+u) \wedge (\theta^c + u^c)
= \omega_0 + (\theta+u) \wedge (\theta^c + u^c)
\]
is the sum of two semi-positive forms. Indeed,
since $u$ is $d^c$-closed, $d^c(\theta+ u) =\omega_0$
is semi-positive; the form  $(\theta+u) \wedge (\theta^c + u^c)$
is semi-positive of rank 1 by definition.

By \ref{_holomo_on_Vaisman_basic_Proposition_}, $u$ is basic.
Since
$\theta+u$ is the sum of $\theta$ and a
basic form, the restriction of
$(\theta+u) \wedge (\theta^c + u^c)$ to $\Sigma$
satisfies
\[
(\theta+u) \wedge (\theta^c + u^c)\restrict \Sigma =
\theta \wedge \theta^c \restrict \Sigma.
\]
The 
sum $\omega_0 + (\theta+u) \wedge (\theta^c + u^c)$
is strictly positive on all tangent
vectors $x\notin \Sigma$ because $\omega_0$ is positive
on these vectors,\footnote{When we say ``a positive
	(1,1)-form $\alpha$ is positive on a vector $v$'',
	we mean that $\1\alpha(v, I(v))>0$; a form is
	Hermitian if it is positive on all non-zero vectors.}
and positive on $x\in \Sigma$
because $\theta \wedge \theta^c \restrict \Sigma$
is positive on such $x$.

\hfill

{\bf Step 2:} It remains to show that
none of the classes $\alpha$ with $\mu(\alpha)\leq 0$ can be
realized as a Lee class of an LCK structure. 
By \ref{_Lee_form_on_Vaisman_is_Vaisman_Proposition_},
any Lee class on $M$ is the Lee class of a Vaisman metric.
If $\mu(\alpha)=0$, we can represent $\alpha$ by
a $d, d^c$-closed form $\alpha_0$. This is impossible by
\ref{_theta_not_d^c_closed_Lemma_}. 

\hfill

{\bf Step 3:} It remains to show that
there are no LCK classes which satisfy
$\mu(\alpha)< 0$. Suppose that 
such a class exists; by \ref{_Lee_form_on_Vaisman_is_Vaisman_Proposition_},
it is the Lee class of a Vaisman manifold, hence
it has an LCK potential. This is impossible,
because two Lee classes for LCK structures with
potential cannot sum to zero, by \ref{_Lee_cannot_be_opposite_Proposition_}.
\endproof

%%%%%%%%%%%%%%%%%%%%%%%%%%%%%%%%%%%%%%%%%%%%%%%%%%%%%%%%%%%%%%%%%%%%%%%%
\subsection{The set of Lee classes on LCK manifolds with potential}
%%%%%%%%%%%%%%%%%%%%%%%%%%%%%%%%%%%%%%%%%%%%%%%%%%%%%%%%%%%%%%%%%%%%%%%%

Now we can prove the main result of this paper.

\hfill

\theorem\label{_Lee_cone_on_LCK-pot_Theorem_}
Let $(M, \theta, \omega)$ be a compact LCK manifold with potential, $\dim_\C M\geq 3$, 
and  $\mu:\; H^1(M, \R)\arrow \R$ a non-zero linear map vanishing
on the space $H^{1,0}(M)\oplus \overline{H^{1,0}(M)}$ which
has codimension 1 by \ref{_LCK_pot_Hodge_decompo_Theorem_}.
Assume that $\mu(\theta) >0$. Then $\xi\in H^1(M, \R)$ is the Lee class of an 
LCK structure with potential on $M$ if and only if $\mu(\xi)>0$.

\hfill

\proof
Let $(M, I_\infty)$ be a Vaisman manifold,
and $\{I_i\}$ the sequence of complex structures
converging to $I_\infty$, such that all manifolds
$(M, I_k)$ are isomorphic to $(M,I)$
(\ref{_Vaisman_limit_of_LCK_pot_Theorem_}).
Given an LCK metric with potential, choose the conformal gauge such that
$\omega_\infty:= d^c\theta_\infty + \theta_\infty \wedge \theta^c_\infty$
on $(M,I_\infty)$ (\ref{_d_theta_c_equation_}). Then the form $ I_kdI_k^{-1}\theta_\infty + \theta_\infty \wedge I_k(\theta_\infty)$
remains strictly positive for almost all manifolds $(M, I_k)$,
because $\lim_k I_k = I_\infty$, and positivity is an open
condition. This implies that $\theta_\infty$ is a Lee form on $(M,I_k)$,
for $k$ sufficiently big.
By \ref{_Lee_cone_on_Vaisman_Theorem_} the set ${\goth L}$ of Lee classes on $(M,I)$ 
contains the half-space $\{u\in H^1(M, \R)\ \ |\ \ \mu_0(u)>0\}$
for some linear map $\mu_0:\; H^1(M, \R)\arrow \R$.
By \ref{_Lee_cannot_be_opposite_Proposition_}, 
${\goth L}$ cannot be bigger than a closed half-space.
However, ${\goth L}$ is open, because the condition
``$d^c\theta + \theta \wedge \theta^c$ is Hermitian'' is open in $\theta$,
hence ${\goth L}$ is an open half-space.
It remains only to show that $\mu_0$ is proportional
to $\mu$. This would follow if we prove that
$\ker \mu=\ker \mu_0$. The space $\ker \mu_0$ is
the set of all classes $\alpha \in H^1(M, \R)$
such that neither $\alpha$ nor $-\alpha$ are Lee classes,
and $\ker \mu$ are classes represented by $d, d^c$-closed forms.
By \ref{_theta_not_d^c_closed_Lemma_},
a Lee class of an LCK manifold cannot be represented
by a $d^c$-closed form, which gives $\ker \mu=\ker \mu_0$.
\endproof

\hfill

\noindent{\bf Acknowledgment:} L.O. thanks Massimiliano
Pontecorvo for very useful discussions during the
preparation of this work, during his visit at
Universit\`a di Roma Tre in June 2021. Both authors thank
Victor Vuletescu for very useful comments on a first
version of the paper.

\hfill

\hfill

{\small
	
	\noindent {\sc Liviu Ornea\\
		University of Bucharest, Faculty of Mathematics and Informatics, \\14
		Academiei str., 70109 Bucharest, Romania}, and:\\
	{\sc Institute of Mathematics ``Simion Stoilow" of the Romanian
		Academy,\\
		21, Calea Grivitei Str.
		010702-Bucharest, Romania\\
		\tt lornea@fmi.unibuc.ro,   liviu.ornea@imar.ro}
	
	\hfill

	\noindent {\sc Misha Verbitsky\\
		{\sc Instituto Nacional de Matem\'atica Pura e
			Aplicada (IMPA) \\ Estrada Dona Castorina, 110\\
			Jardim Bot\^anico, CEP 22460-320\\
			Rio de Janeiro, RJ - Brasil }\\
		also:\\
		Laboratory of Algebraic Geometry, \\
		Faculty of Mathematics, National Research University 
		Higher School of Economics,
		6 Usacheva Str. Moscow, Russia}\\
	\tt verbit@verbit.ru, verbit@impa.br }


{\small

\begin{thebibliography}{100}
	
\bibitem[AS]{andreotti_siu} A. Andreotti, Y.T. Siu,  {\em Projective embeddings of pseudoconcave spaces}, Ann. Scuola Norm. Sup. Pisa {\bf 24}, 231--278 (1970).
	
\bibitem[AD]{ad2} V. Apostolov, G. Dloussky, {\em On the Lee classes of locally conformally symplectic complex surfaces}, J. Sympl. Geom.  {\bf 16} (2018), 931--958. arXiv:1611.00074.

\bibitem[Be]{bel} F.A. Belgun, {\em On the metric structure of non-K\"ahler complex surfaces}, Math. Ann. {\bf 317} (2000), 1--40.

\bibitem[BT]{_Borel_Tits:Groupes_Reductifs_} A. Borel, J. Tits, {\em Groupes r\'eductives}, 
Publications Math\'ematiques de l'IH\'ES, {\bf 27} (1965),  55--151.



\bibitem[Br]{_Brunella:Kato_} M. Brunella, {\em Locally conformally K\"ahler metrics on Kato surfaces}, Nagoya Math. J. {\bf 202} (2011), 77--81.

\bibitem[D]{_Demailly:Book_} 
J.-P. Demailly, Complex analytic and differential geometry, {\small\url{ https://www-fourier.ujf-grenoble.fr/~demailly/manuscripts/agbook.pdf}}


\bibitem[DO]{do} S. Dragomir, L. Ornea, Locally conformally K\"ahler manifolds, Progress in Math. {\bf 55}, Birkh\"auser, 1998.

\bibitem[EM]{em} Y. Eliashberg, E. Murphy, {\em Making cobordisms symplectic}, arXiv:1504.06312.

\bibitem[F]{friedman} A. Friedman, Foundations of modern analysis, Dover, 2010.

\bibitem[GO]{go} P. Gauduchon, L. Ornea, {\em Locally conformally K\"ahler metrics on Hopf surfaces}, Ann. Inst. Fourier (Grenoble) {\bf 48} (1998), 1107--1128.

\bibitem[Hum]{hum} J.E. Humphreys, Linear algebraic groups, GTM 21, 4th ed., Springer, 1998.

\bibitem[I]{_Istrati:LCK-pot_} N. Istrati, {\em Existence
  criteria for special locally conformally K\"ahler
  metrics}, Ann. Mat. Pura. Appl. {\bf 198} (2019),
  335--353.


\bibitem[IOP]{iop}
N. Istrati, A. Otiman, M. 
Pontecorvo, {\em On a class of Kato manifolds},
  IMRN, {\bf 7} (2021), 5366--5412. arXiv:1905.03224.

\bibitem[IOPR]{iopr} N. Istrati, A. Otiman, M. 
Pontecorvo, M. Ruggiero, {\em Toric Kato manifolds}, arXiv:2010.14854.

\bibitem[KO]{kor} Y. Kamishima, L. Ornea, {\em Geometric flow on compact  locally conformally K\"ahler manifolds}, Tohoku Math. J., {\bf 57} (2) (2005), 201--221.

\bibitem[Ka]{kashiwada_kodai} T. Kashiwada, {\em On V-harmonic forms in compact locally conformal K\"ahler manifolds with the parallel Lee form},
Kodai Math. J. {\bf 3} (1980) 70--82.

%\bibitem[MaMoPi]{mamopi} F. Madani, A. Moroianu, M. Pilca, {\em On toric locally conformally K\"ahler manifolds}, Ann. Global. Anal. Geom.,  {\bf 51} (2017), no. 4, 401--417.

\bibitem[OT]{ot} K. Oeljeklaus, M. Toma, {\em Non-K\"ahler compact complex manifolds associated to number  fields}, Ann. Inst. Fourier {\bf 55}, no. 1 (2005), 1291--1300.

%\bibitem[Ok]{_Oka_} K. Oka, {\em Domaines finis sans point critique int\'erieur}, Japanese J. Math. {\bf 23} (1953), 97--155.
	
\bibitem[OV1]{ov_jgp_09} 
L. Ornea, M. Verbitsky, {\em
  Morse-Novikov cohomology of locally conformally K\"ahler
  manifolds}, J. Geom. Phys. {\bf 59}, No. 3 (2009),
  295--305.


\bibitem[OV2]{ov_lckpot} 
L. Ornea, M. Verbitsky, {\em
  Locally conformal K\"ahler manifolds with potential},
  Math. Ann. {\bf 348} (2010), 25--33.


\bibitem[OV3]{ov_imrn_10} 
L. Ornea, M. Verbitsky, {\em Topology of Locally
  Conformally K\"ahler Manifolds with Potential}, IMRN,
Vol. 2010, pp. 717--726.

\bibitem[OV4]{ov_jgp_16} 
L. Ornea, M. Verbitsky, {\em LCK rank of locally conformally K\"ahler manifolds with potential}, J. Geom. Phys. {\bf 107} (2016), 92--98.

\bibitem[OV5]{ov_pams} 
L. Ornea, M. Verbitsky, {\em Locally conformally K\"ahler
  metrics obtained from pseudoconvex shells},
Proc. Amer. Math. Soc. {\bf 144} (2016), 325--335.

	
\bibitem[OV6]{_ov_super_sas_}  L. Ornea, M. Verbitsky,
  {\em Supersymmetry and Hodge theory on Sasakian and
    Vaisman manifolds}, arXiv:1910.01621v2. To appear in Manuscripta Math.

\bibitem[Ot]{oti2} A. Otiman, {\em Morse-Novikov cohomology of locally conformally K\"ahler surfaces},   Math. Z. {\bf 289} (2018), no. 1-2, 605--628. arXiv:1609.07675.

\bibitem[Ros]{rossi} H. Rossi, {\em Attaching analytic spaces to an analytic space along a pseudo--convex boundary}, Proceedings of the Conference Complex Manifolds (Minneapolis), pp. 242--256. Springer, Berlin (1965).

\bibitem[Ts1]{tsuk} K. Tsukada, {\em Holomorphic forms and holomorphic vector fields on compact generalized Hopf manifolds}, Compositio Math. {\bf 93} (1994), no. 1, 1--22.

\bibitem[Ts2]{tsu} K. Tsukada, {\em The canonical foliation of a compact generalized Hopf manifold}, Differential Geom. Appl. {\bf 11} (1999), no. 1, 13--28.

\bibitem[Va1]{va_tr} I. Vaisman, {\em On locally and globally conformal K\"ahler manifolds}, Trans. Amer. Math. Soc., {\bf 262} (1980), 533-542.
	
\bibitem[Va2]{va_gd} 
I. Vaisman, {\em Generalized Hopf manifolds}, Geom. Dedicata, {\bf 13} (1982), 231--255.

\bibitem[Ve]{_Verbitsky:Vanishing_LCHK_} 
M. Verbitsky, {\em Theorems on the vanishing of 
cohomology for locally conformally hyper-K\"ahler
manifolds},  Proc. Steklov Inst. Math. no. 3 ({\bf 246}), 54--78 (2004).


\bibitem[VVO]{_ovv:surf_} 
M. Verbitsky, V. Vuletescu, L. Ornea {\em Classification of non-K\"ahler surfaces and 
locally conformally K\"ahler geometry}, 
Russian Math. Surv. {\bf 76} (2021),
261--290. arxiv:1810.05768.

\bibitem[Vu]{_Vuletescu:blowups_} V. Vuletescu, {\em Blowing-up points on l.c.K. manifolds}, Bull. Math. Soc. Sci. Math. Roumanie (N.S.) {\bf 52(100)} (2009), no. 3, 387--390

\end{thebibliography}
}
\end{document}